\documentclass[a4paper, 11pt]{amsart}

\allowdisplaybreaks

\usepackage{cite}
\usepackage[mathscr]{eucal}
\usepackage{amsmath}
\usepackage{amssymb}
\usepackage{amsfonts}
\usepackage{amsthm}

\usepackage{cases}

\usepackage[mathscr]{eucal}
\usepackage{amsmath}
\usepackage{amssymb}
\usepackage{amsfonts}
\usepackage{amsthm}
\usepackage[dvips]{graphics, color}

\theoremstyle{plain}
\newtheorem{theorem}{Theorem}[section]
\newtheorem{proposition}{Proposition}[section]

\theoremstyle{remark}

\newtheorem{remark}{Remark}[section]

\numberwithin{equation}{section}

\def\<{\left<} \def\>{\right>}

\def\proof{\noindent{\it Proof. }}
\def\bea{\begin{eqnarray} }
\def\eea{\end{eqnarray} }
\def\be{\begin{equation} }
\def\ee{\end{equation} }

\begin{document}

\title[Hamiltonian stationary Lagrangian surfaces]{Hamiltonian stationary Lagrangian surfaces with  harmonic  mean curvature  in complex space forms}
\author[T. Sasahara]{Toru Sasahara}
\address{Division of Mathematics, Center for Liberal Arts and Sciences, 
Hachinohe Institute of Technology, 
Hachinohe, Aomori, 031-8501, Japan}
\email{sasahara@hi-tech.ac.jp}

\date{}

\begin{abstract}
In this paper, we study Hamiltonian stationary Lagrangian surfaces in complex space forms. We first show that when the mean curvature is a non-zero constant, the second fundamental form is parallel.
We then consider the case in which the mean curvature is a non-constant harmonic function. Under the additional assumption that the Gaussian curvature is constant, we obtain a complete classification of such Lagrangian surfaces.
\end{abstract}


 \maketitle

 \section{Introduction}
 Let $\tilde M^n$ be a
complex $n$-dimensional K\"ahler manifold with  complex structure $J$ 
and K\"ahler metric $\<\,, \>$.
 An $n$-dimensional submanifold $M$ of  $\tilde M^n$ is called 
 Lagrangian if  $\<X, JY\>=0$ for all tangent vector fields $X$, $Y$ on $M$.
A normal vector field $\xi$ of a Lagrangian submanifold $M$  is called a {Hamiltonian variation} if
 $\xi=J\nabla f$ for some compactly supported  function $f$ on $M$, 
 where $\nabla$ is the gradient on $M$.
 A Lagrangian submanifold   is said to be {Hamiltonian stationary}  if 
 it is a critical point of the volume functional for all Hamiltonian variations.
As shown by Oh  \cite{oh}, this condition is equivalent to  the divergence-free condition
   \be {\rm div}(JH)=0,\label{HS}\ee 
 where $H$ denotes the mean curvature vector field of $M$.

In this paper, we investigate  Hamiltonian stationary Lagrangian surfaces in complex space forms. When the mean curvature is  a non-zero constant, we show that the second fundamental form is parallel. Such Lagrangian surfaces have  been completely classified in  Theorem 7.2 of \cite{c2}, Theorems A.1 and A.2 of \cite{cdv}.
We then turn to the case in which the mean curvature is a non-constant harmonic function. 
Assuming that the Gaussian curvature is constant, 
we prove that the ambient space must be  the complex hyperbolic plane and 
the Gaussian curvature necessarily takes a specific negative value.
This phenomenon leads to a complete classification.


\section{Preliminaries}

  Let $\tilde M^n(4\epsilon)$ be a complete and simply connected complex space form of complex dimension $n$ and constant
 holomorphic sectional curvature $4\epsilon$, that is, 
  $\tilde M^n(4\epsilon)$ is  
 the complex Euclidean space $\mathbb{C}^n$, 
 the complex projective space $\mathbb{C}P^n(4\epsilon)$ or the complex hyperbolic space $\mathbb{C}H^n(4\epsilon)$
 according as  $\epsilon=0$, $\epsilon>0$ or $\epsilon<0$.

 Let $M$ be   a Lagrangian submanifold of  $\tilde M^n(4\epsilon)$.
 We denote the Levi-Civita connections on $M^n$ and 
$\tilde M^n(4\epsilon)$ by $\nabla$ and
 $\tilde\nabla$, respectively. The
 Gauss and Weingarten formulas are  given respectively by
\be
 \tilde \nabla_XY = \nabla_XY+h(X,Y), \quad 
 \tilde\nabla_X \xi = -S_{\xi}X+D_X\xi \nonumber
\ee
 for tangent vector fields $X$, $Y$ and normal vector field $\xi$,
 where $h$, $S$ and $D$ are the second fundamental 
 form, the shape operator and the normal
 connection. 
  The mean curvature vector field $H$ is defined by 
$H=(1/n){\rm trace}\,h.$
 The function $|H|$ is called the  mean curvature.
 We have (cf. \cite{co})
 \begin{align}
 &  D_XJY=J(\nabla_XY), \label{DX}\\
 &  \<h(X, Y), JZ\>=\<h(Y, Z), JX\>=\<h(Z, X), JY\>. \label{AX}
 \end{align}

 Denote by $R$  the Riemann curvature tensor of $\nabla$.
 Then the equations of
 Gauss and  Codazzi are given  respectively by
 \begin{align}
 \<R(X,Y)Z,W\>=&
 \<h(Y, Z), h(X, W)\>-\<h(X, Z), h(Y, W)\>\nonumber \\
 &+\epsilon(\<Y,Z\>\<X,W\>-\<X,Z\>\<Y,W\>), \label{Gauss}\\
 ({\bar\nabla}_{X}h)(Y,Z)=&
 ({\bar\nabla}_{Y}h)(X,Z),\label{Codazzi}
 \end{align}
 where $X,Y,Z,W$ 
  are vectors tangent  to
 $M$,  and
 $\bar\nabla h$ is defined by
 \be ({\bar\nabla}_{X}h)(Y,Z)= D_X h(Y,Z) - h(\nabla_X
 Y,Z) - h(Y,\nabla_X Z). \nonumber
 \ee

 \section{Hamiltonian stationary Lagrangian surfaces with harmonic mean curvature}
 Let $M$  be a  Lagrangian surface 
in $\tilde M^2(4\epsilon)$, where $\epsilon\in\{-1, 0, 1\}$.
 Suppose that $H\ne 0$ everywhere. Denote by  $K$ the Gaussian curvature of $M$.
Let $\{e_1, e_2\}$ be a local  orthonormal frame on $M$ such that 
$Je_1$ is parallel to $H$. 
It follows from (\ref{AX}) that the
 second fundamental form takes the form
\be
\begin{split} \label{sf}
&  h(e_1, e_1)=(a-c)Je_1+bJe_2,\\
&  h(e_1, e_2)=bJe_1+cJe_2, \\
&  h(e_2, e_2)=cJe_1-bJe_2 
\end{split}
\ee
for some functions $a$, $b$ and  $c$.

Putting $\omega_i^j(e_k)=\<\nabla_{e_k}e_i, e_j\>$, and using (\ref{DX}) and (\ref{sf}),  we have
\begin{align}
(\bar\nabla_{e_1}h)(e_2, e_2)&=\{e_1c+3b\omega_1^2(e_1)\}Je_1
-\{e_1b-3c\omega_1^2(e_1)\}Je_2,\nonumber\\
(\bar\nabla_{e_2}h)(e_1, e_2)&=\{e_2b+(a-3c)\omega_1^2(e_2)\}Je_1
+\{e_2c+3b\omega_1^2(e_2)\}Je_2,\nonumber\\
(\bar\nabla_{e_1}h)(e_1, e_2)&=\{e_1b+(a-3c)\omega_1^2(e_1)\}Je_1
+\{e_1c+3b\omega_1^2(e_1)\}Je_2,\nonumber\\
(\bar\nabla_{e_2}h)(e_1, e_1)&=\{e_2(a-c)-3b\omega_1^2(e_2)\}Je_1
+\{e_2b+(a-3c)\omega_1^2(e_2)\}Je_2. \nonumber
\end{align}
Therefore, 
 the Codazzi equation (\ref{Codazzi})  implies 
\begin{align}
& e_1c+3b\omega_1^2(e_1)=e_2b+(a-3c)\omega^2_1(e_2), \label{C1}\\
& -e_1b+3c\omega_1^2(e_1)=e_2c+3b\omega^2_1(e_2),\label{C2}\\
& e_2(a-c)-3b\omega_1^2(e_2)=e_1b+(a-3c)\omega^2_1(e_1).\label{C3}
\end{align}
Combining (\ref{C2}) and (\ref{C3}) yields 
\be
e_2a-a\omega^2_1(e_1)=0.\label{C4}
\ee

Assume that $M$ is Hamiltonian stationary. Then, by (\ref{HS}) we have   
\be e_1a+a\omega_1^2(e_2)=0. \label{C5}
\ee
Using (\ref{C4}) and (\ref{C5}), we obtain
\be [a^{-1}e_1, a^{-1}e_2]=0. \nonumber
\ee
Therefore, there exists a  local coordinate system $\{u, v\}$ such that
$ e_1=a\partial_u$ and $e_2=a\partial_v$.
  Hence, the metric tensor is given by 
\be 
g=a^{-2}(du^2+dv^2), \label{g}
\ee 
which implies that 
\begin{align}
&\omega_1^2(e_1)=a_v, \quad 
\omega_1^2(e_2)=-a_u,\label{omega} \\
&K=-(a_u)^2-(a_v)^2+a(a_{uu}+a_{vv}). \label{K}
\end{align}
It follows from (\ref{omega}) that (\ref{C1}) and (\ref{C2}) can be rewritten as 
\begin{align}
ac_u+3ba_v&=ab_v-(a-3c)a_u, \label{au}\\
  -ab_u+3ca_v&=ac_v-3ba_u.\label{av}
\end{align}
Put $\delta=\epsilon-K$. Then the Gauss equation  (\ref{Gauss}) together with  (\ref{sf}) yields
\be
\delta=2b^2-ac+2c^2. \label{G1}
\ee

 In the case where the mean curvature is constant, we have

 \begin{proposition}\label{prop}
 Let $M$ be a Hamiltonian stationary Lagrangian surface 
in $\tilde M^2(4\epsilon)$.
If $|H|$ is a non-zero constant, then $\bar\nabla h=0$.
 \end{proposition}
 \proof
 From (\ref{omega}), it follows that  if $|H|=|a|/2$  is a non-zero constant, then  $\omega_1^2=0$, and hence $K=0$. 
 Equations (\ref{C1}) and (\ref{C2}) reduce to 
\be 
c_u-b_v=c_v+b_u=0.\label{CR}
\ee
Since $\delta=\epsilon$, differentiating  (\ref{G1}) with respect to $u$ and $v$,  and using  (\ref{CR}), we have
\be
    \begin{pmatrix}
       4b & 4c-a \\
       a-4c & 4b\\
   \end{pmatrix}
  \begin{pmatrix}
  b_u \\
  c_u
  \end{pmatrix}=\begin{pmatrix}
  0 \\
  0
  \end{pmatrix}.\label{mt}
\ee

If $16b^2+(a-4c)^2\ne 0$ on an open subset, then by (\ref{mt}) we obtain $b_u=c_u=0$. Thus, it follows from (\ref{CR}) that
$b$ and $c$ are constant. Otherwise, $b=0$ and $c$ is constant.
In both cases, using  (\ref{DX}) and (\ref{sf}), we obtain $\bar\nabla h=0$. \qed

\begin{remark}
    An explicit description of all Lagrangian surfaces with $\bar\nabla h=0$ in  $\tilde M^2(4\epsilon)$
  has been obtained in Theorem 7.2 of \cite{c2}, Theorems A.1 and A.2 of \cite{cdv}. 
\end{remark}


Our next step is to investigate the case of non-constant mean curvature.
The main result of this paper is the following.

\begin{theorem}\label{main} 
Let $M$ be a Hamiltonian stationary Lagrangian surface 
in a complex space form $\tilde M^2(4\epsilon)$, where $\epsilon\in\{-1, 0, 1\}$.
Suppose that $H$ is nowhere vanishing. 
If $|H|$ is a  non-constant harmonic function and $K$ is constant, then $K=\epsilon=-1$ and  $M$ is locally congruent to
the image of  $\Pi\circ \phi$, where 
$\Pi: H_1^5(-1)\rightarrow {\mathbb C}H^2(-4)$
is the Hopf fibration
and 
 $\phi: M\rightarrow H_1^5(-1)\subset {\mathbb C}^3_1$   is given by one of the following immersions$:$

 $(1)$ 
\begin{align*}
\phi (x, y)=\biggl(m e^{y}+\dfrac{e^{-y}+2im^2xe^{y}}{2m}, m e^{ix+y},
\dfrac{e^{-y}+2im^2xe^{y}}{2m}\biggr);
\end{align*}

$(2)$ 
\be 
\phi (x,y)=\biggl(1-\dfrac{i(1+m^2)}{m^2x+y}, \dfrac{m\sqrt{1+m^2}e^{ix}}{m^2x+y}, 
\dfrac{\sqrt{1+m^2}e^{iy}}{m^2x+y}\biggr),\nonumber
\ee
where 
 $m$ is a positive real number. 
\end{theorem}

\proof 
Suppose that $|H|$ is a harmonic function on $M$. Then by  (\ref{g}) we have
\be a_{uu}+a_{vv}=0.\label{har}
\ee
Moreover, suppose that  $|H|$ is non-constant and $K$ is constant.  Then, combining  (\ref{K}) and (\ref{har}) shows that  $K<0$ and 
\be
 a_u=\sqrt{-K}\cos\theta, \quad 
 a_v=\sqrt{-K}\sin\theta.  \label{auv}
\ee
for some function $\theta$ on $M$. 
Substituting (\ref{auv}) into (\ref{har}), we get
\be -(\sin\theta)\theta_u+(\cos\theta)\theta_v=0. \label{theta1}\ee

On the other hand, since $a_{uv}-a_{vu}=0$ holds, by (\ref{auv}) we obtain
\be (\cos\theta)\theta_u+(\sin\theta)\theta_v=0.\label{theta2}\ee 
It follows from (\ref{theta1}) and  (\ref{theta2}) that $\theta_u=\theta_v=0$, that is, $\theta$ is constant.
Solving (\ref{auv}), we conclude that up to translations, $a$ is given by
\be
a=\sqrt{-K}\{(\cos\theta)u+(\sin\theta)v\}. \label{a}
\ee

{\bf Case (i):} $b=0$ on an open subset $\mathcal{U}$. In this case,   
 (\ref{au}) and (\ref{av}) reduce respectively to
 \begin{align}
 & (a-3c)a_u+ac_u=0, \label{C6}\\
 &  3ca_v-ac_v=0.\label{C7}
 \end{align}
Differentiating (\ref{G1}) with respect to $v$, we have
\be
ca_v+(a-4c)c_v=0.\label{D1}
\ee
Combining (\ref{C7}) and  (\ref{D1}) gives 
\be 
(a-c)c_v=0.\nonumber
\ee

If $c_v\ne 0$ on an open subset in $\mathcal{U}$, then $a=c$. From
(\ref{G1}) we see that $a$ is constant, which contradicts our assumption.  
Hence we have $c_v=0$ on $\mathcal{U}$. Thus, by (\ref{C7}) we get 
$$ca_v=0.$$

{\bf Case (i.1):} $c=0$ on an open subset $\mathcal{U}_1\subset\mathcal{U}$. In this case,
from (\ref{G1}) we have  $K=\epsilon=-1$.
It follows from  (\ref{C6}) that  $a_u=0$ on $\mathcal{U}_1$. Hence, using (\ref{a}) and $K=-1$, we obtain $a^2=v^2$.  Applying  the coordinate transformation $y=-\int v^{-1}dv$,
we see that the metric tensor (\ref{g}) becomes
\be
g=m^2e^{2y}du^2+dy^2\nonumber
\ee
for some positive constant $m$, and the second fundamental form satisfies
\be  
h(\partial_u, \partial_u)=J\partial_u, \quad h(\partial_u, \partial_y)=h(\partial_y,\partial_y)=0.\nonumber
\ee
We rewrite $u$ as $x$. According to \cite[page 3475]{cd},  we conclude that $\mathcal{U}_1$ is  congruent to the Lagrangian surface obtained from (1).

{\bf Case (i.2):} $a_v=0$ on an open subset $\mathcal{U}_2\subset\mathcal{U}$. In this case, $a_u\ne 0$.
Differentiating (\ref{G1}) with respect to $u$ leads to
\be
ca_u+(a-4c)c_u=0.\label{D2}
\ee
Eliminating $c_u$ from  (\ref{C6}) and (\ref{D2}) gives 
\be
(a-2c)(a-6c)a_u=0.\nonumber
\ee

If $a\ne 2c$ on an open subset of $\mathcal{U}_2$, then $a=6c$, and it follows from (\ref{G1}) that $a$ is constant,
which contradicts our assumption. Hence $a=2c$, which together with (\ref{G1}) implies  $K=\epsilon=-1$.  Thus,  by (\ref{a}) we obtain 
$a^2=u^2$.
Applying  the coordinate transformation $x=(u+v)/2$ and $y=(u-v)/2$, we see that
the metric tensor (\ref{g}) becomes
\be
g=\dfrac{2}{(x+y)^2}(dx^2+dy^2)\nonumber
\ee
 and the second fundamental form satisfies
\be  
h(\partial_x, \partial_x)=J\partial_x, \quad h(\partial_x, \partial_y)=0, \quad
h(\partial_y, \partial_y)=J\partial_y.\nonumber
\ee
According to \cite[page 124]{cdvv}, $\mathcal{U}_2$ is congruent to
the Lagrangian surface obtained from (2) with $m=1$.

{\bf Case (ii):} $b\ne 0$ on an open subset $\mathcal{V}$.
We put $A=a_u$ and $B=a_v$, which are constant.
Differentiating (\ref{G1}) with respect to $u$ and $v$, we obtain
\begin{align}
& b_u=\frac{cA+(a-4c)c_u}{4b}, \label{bu}\\
&  b_v=\frac{cB+(a-4c)c_v}{4b}.\label{bv}
\end{align}
Substituting (\ref{bu}) and (\ref{bv}) into (\ref{au}) and (\ref{av}) yields
\begin{align}
&  a\{4bc_u-(a-4c)c_v\}=4b(3c-a)A-(12b^2-ac)B, \label{cu}\\
& a\{(a-4c)c_u+4bc_v\}=(12b^2-ac)A+12bcB. \label{cv}
\end{align}
Solving (\ref{cu}) and (\ref{cv}) for $c_u$ and $c_v$, we have
\begin{align}
& c_u=\frac{Af_1+Bf_2}{a(a^2+16b^2-8ac+16c^2)},\label{Cu}\\
&  c_v=\frac{Af_3+Bf_4}{a(a^2+16b^2-8ac+16c^2)},\label{Cv}
\end{align}
where
\begin{align}
& f_1=-4ab^2-a^2c+4ac^2,  \\
& f_2=-48b^3+16abc-48bc^2,\\
& f_3=4a^2b+48b^3-32abc+48bc^2,  \\
& f_4=12ab^2-a^2c+4ac^2.
\end{align}
Substituting (\ref{Cu}) and (\ref{Cv}) into (\ref{bu}) and (\ref{bv}) gives
\begin{align}
b_u=& \frac{Ag_1+Bg_2}{a(a^2+16b^2-8ac+16c^2)},\label{Bu}\\
b_v=&\frac{Ag_3+Bg_4}{a(a^2+16b^2-8ac+16c^2)},\label{Bv}
\end{align}
where 
\begin{align*}
& g_1=-a^2 b + 8 a b c, \\
& g_2=-12ab^2+4a^2c+48b^2c
-28ac^2+48c^3,\\
&g_3=a^3+12ab^2-12a^2c-48b^2c+44ac^2-48c^3,\\
&g_4=3 a^2 b - 8 a b c.
\end{align*}

Using  (\ref{Cu})-(\ref{Bv}), we differentiate (\ref{Cu}) and (\ref{Cv}) with respect to $v$ and $u$, respectively.  
A computation with a computer algebra system yields
\begin{align}
&c_{uv}=\frac{A^2h_1+ABh_2+B^2h_3}{a^2(a^2+16b^2-8ac+16c^2)^2}, \label{cuv}\\
&c_{vu}=\frac{A^2h_4+ABh_5+B^2h_6}{a^2(a^2+16b^2-8ac+16c^2)^2}, \label{cvu}
\end{align}
where 
\begin{align*}
h_1=&-12 a^4 b - 144 a^2 b^3 + 160 a^3 b c + 768 a b^3 c - 656 a^2 b c^2 + 
 768 a b c^3,\\
   h_2=&-108 a^3 b^2 - 960 a b^4 + 18 a^4 c + 1168 a^2 b^2 c + 2304 b^4 c - 
 252 a^3 c^2\\
 &-3840 a b^2 c^2 + 1296 a^2 c^3 + 4608 b^2 c^3 - 2880 a c^4 + 2304 c^5,\\
    h_3=&-96 a^2 b^3 + 768 b^5 - 384 a b^3 c + 32 a^2 b c^2 + 1536 b^3 c^2 - 
 384 a b c^3 + 768 b c^4,\\
    h_4=&-8 a^4 b - 96 a^2 b^3 - 768 b^5 + 128 a^3 b c + 1152 a b^3 c - 
 608 a^2 b c^2\\
   &-1536 b^3 c^2 + 1152 a b c^3 - 768 b c^4,\\
   h_5=&-92 a^3 b^2 - 192 a b^4 + 18 a^4 c + 784 a^2 b^2 c + 2304 b^4 c - 
 252 a^3 c^2\\
   &-3072 a b^2 c^2 + 1296 a^2 c^3 + 4608 b^2 c^3 - 2880 a c^4 + 2304 c^5,\\
   h_6=&-240 a^2 b^3 + 80 a^3 b c + 768 a b^3 c - 496 a^2 b c^2 + 768 a b c^3.
\end{align*}
Subtracting    (\ref{cvu}) from (\ref{cuv}), we have 
\be  
c_{uv}-c_{vu}=\frac{4b(A^2k_1+ABk_2+B^2k_3)}{a^2(a^2+16b^2-8ac+16c^2)^2},\label{integ1}
\ee
where
\begin{align*}
k_1=&-a^4 - 12 a^2 b^2 + 192 b^4 + 8 a^3 c - 96 a b^2 c - 12 a^2 c^2\\
& + 384 b^2 c^2 - 96 a c^3 + 192 c^4,\\
k_2=&-4 a^3 b - 192 a b^3+ 96 a^2 b c - 192 a b c^2,\\
k_3=&36 a^2 b^2 + 192 b^4 - 20 a^3 c - 288 a b^2 c\\
 &+ 132 a^2 c^2 + 
   384 b^2 c^2 - 288 a c^3 + 192 c^4.
\end{align*}
Since the compatibility condition  $c_{uv}=c_{vu}$ must hold, it follows from $(\ref{integ1})$
 and the assumption $b\ne 0$ that
\be 
A^2k_1+ABk_2+B^2k_3=0. \label{integ2}
\ee
Differentiating the left hand side of (\ref{integ2}) and using (\ref{Cu})-(\ref{Bv}), with the aid of
a computer algebra system, we get
\begin{align} 
(A^2k_1+ABk_2+B^2k_3)_u&=\frac{4(A^3P_1+A^2BP_2+AB^2P_3+B^3P_4)}{a^2+16b^2-8ac+16c^2},\label{integ2u}\\
(A^2k_1+ABk_2+B^2k_3)_v&=\frac{8(A^3P_5+A^2BP_6+AB^2P_7+B^3P_8)}{a^2+16b^2-8ac+16c^2},\label{integ2v}
\end{align}
where
\begin{align*}
P_1=&-a^5 - 24 a^3 b^2 - 192 a b^4 + 12 a^4 c + 168 a^2 b^2 c+  384 b^4 c \\
& - 56 a^3 c^2 - 576 a b^2 c^2 + 168 a^2 c^3 + 
   768 b^2 c^3 - 384 a c^4 + 384 c^5,\\
   P_2=&-2 a^4 b - 72 a^2 b^3 - 1920 b^5 + 24 a^3 b c+ 1344 a b^3 c\\
   & - 264 a^2 b c^2
    - 3840 b^3 c^2 + 1344 a b c^3 - 1920 b c^4,\\
   P_3=&32 a^3 b^2 + 960 a b^4 - 14 a^4 c - 888 a^2 b^2 c - 
   2688 b^4 c + 224 a^3 c^2\\
    &+ 4032 a b^2 c^2 - 1272 a^2 c^3 - 
   5376 b^2 c^3 + 3072 a c^4 - 2688 c^5,\\
   P_4=&24 a^2 b^3 + 1152 b^5 - 8 a^3 b c - 960 a b^3 c \\
   &+216 a^2 b c^2 + 2304 b^3 c^2 - 960 a b c^3 + 1152 b c^4,\\
   P_5=&a^4 b + 60 a^2 b^3 + 576 b^5 - 32 a^3 b c - 672 a b^3 c \\
&+ 252 a^2 b c^2 + 1152 b^3 c^2 - 672 a b c^3 + 576 b c^4,\\
   P_6=&-a^5 - 38 a^3 b^2 - 192 a b^4+ 24 a^4 c + 444 a^2 b^2 c
   + 1344 b^4 c - 218 a^3 c^2 \\
   &- 2016 a b^2 c^2 + 924 a^2 c^3 + 
   2688 b^2 c^3 
   - 1824 a c^4 + 1344 c^5,\\
   P_7=&-4 a^4 b - 180 a^2 b^3 - 960 b^5 + 96 a^3 b c+ 1248 a b^3 c\\
   & - 
   564 a^2 b c^2 - 1920 b^3 c^2 + 1248 a b c^3 - 960 b c^4,\\
   P_8=&6 a^3 b^2 - 5 a^4 c - 84 a^2 b^2 c - 192 b^4 c + 50 a^3 c^2\\
   & + 288 a b^2 c^2 - 180 a^2 c^3 - 384 b^2 c^3 + 288 a c^4 - 192 c^5.
\end{align*}
Thus,  it follows from (\ref{integ2}), (\ref{integ2u}) and (\ref{integ2v})   that 
\begin{align}
& A^3P_1+A^2BP_2+AB^2P_3+B^3P_4=0, \label{integ3}\\
& A^3P_5+A^2BP_6+AB^2P_7+B^3P_8=0. \label{integ4}
\end{align}

Using a computer algebra system, we find that
the resultant of the left-hand sides of (\ref{integ3}) and (\ref{integ4}) with respect to $A$ is, up to a non-zero constant factor,
\be
B^9(2b^2-ac+2c^2)^3(a^2 + 16 b^2 - 8 a c + 16 c^2)^7
 (a^2 + 48 b^2 - 24 a c + 
   48 c^2)^2Q,\label{R1}
\ee
where
\begin{align*}
 Q=&5 a^6 + 120 a^4 b^2 + 720 a^2 b^4 + 2304 b^6 - 
   60 a^5 c - 720 a^3 b^2 c \\
   &- 3456 a b^4 c + 300 a^4 c^2 + 
   3168 a^2 b^2 c^2 + 6912 b^4 c^2 - 1008 a^3 c^3 \\
   &- 6912 a b^2 c^3 + 
   2448 a^2 c^4 + 6912 b^2 c^4 - 3456 a c^5 + 2304 c^6.
\end{align*}
Combining (\ref{G1}) and  (\ref{R1}), 
we find that (\ref{R1}) simplifies to
\be 
B^9 \delta^3 (a^2 + 8 \delta)^7 (a^2 + 24 \delta)^2 (5 a^6 + 60 a^4 \delta + 
   180 a^2 \delta^2 + 288 \delta^3),\label{R2}\nonumber
\ee
which must vanish on $\mathcal{V}$. Taking  into account that $\nabla a\ne 0$, we have
\be B\delta=0.\nonumber
\ee

If $\delta\ne 0$, then  $B=0$, which implies that
(\ref{integ2}) and (\ref{integ4}) reduce respectively to
\begin{align}
&a^4 + 12 a^2 b^2 - 192 b^4 - 8 a^3 c + 96 a b^2 c 
+ 12 a^2 c^2 \nonumber\\
&- 
   384 b^2 c^2 + 96 a c^3 - 192 c^4=0,\label{abc}\\
  &a^4 + 60 a^2 b^2 + 576 b^4 - 32 a^3 c - 672 a b^2 c \nonumber \\
&+252 a^2 c^2 + 1152 b^2 c^2 - 672 a c^3 + 576 c^4=0.\label{abc2}
\end{align}
Subtracting   (\ref{abc}) from (\ref{abc2}), we derive
\be
24 (2 b^2 - a c + 2 c^2) (a^2 + 16 b^2 - 8 a c + 16 c^2)=0. \label{abc3}
\ee
Combining (\ref{G1}) and  (\ref{abc3}) implies
\be
a^2+8\delta=0,\nonumber
\ee
which contradicts $\nabla a\ne 0$. Hence  $\delta=0$, that is, $K=\epsilon=-1$.

From (\ref{G1}) and $b\ne 0$, we see that $c\ne0$. Changing the sign of $e_1$ if necessary, we may assume that $c>0$. Put
\begin{align}
f=\frac{\sqrt{b^2+c^2-b\sqrt{b^2+c^2}}}{\sqrt{2}(b^2+c^2)},\quad 
k=\frac{\sqrt{b^2+c^2+b\sqrt{b^2+c^2}}}{\sqrt{2}(b^2+c^2)}, \label{fk}
\end{align}
which are non-zero  and unequal everywhere.
It follows from   (\ref{G1}) and (\ref{fk}) that
\begin{align}
a=\frac{\sqrt{f^2+k^2}}{fk},\quad 
b=\frac{k^2-f^2}{(f^2+k^2)^{\frac{3}{2}}},\quad
c=\frac{2fk}{(f^2+k^2)^{\frac{3}{2}}}.\label{abcfk}
\end{align}
We make the following change of  basis:
\begin{align}
\tilde e_1=\frac{k}{\sqrt{f^2+k^2}}e_1+\frac{f}{\sqrt{f^2+k^2}}e_2, \quad
\tilde e_2=\frac{f}{\sqrt{f^2+k^2}}e_1-\frac{k}{\sqrt{f^2+k^2}}e_2. \label{change}
\end{align}
Then, using (\ref{sf}), (\ref{abcfk}) and   (\ref{change}),  we find
\begin{align*}
h(\tilde e_1, \tilde e_1)=f^{-1}J\tilde e_1,\quad 
 h(\tilde e_1, \tilde e_2)=0, \quad
  h(\tilde e_2, \tilde e_2)=k^{-1}J\tilde e_2.
\end{align*}

According to \cite[Theorem 4.1]{cdvv},  there exists a local coordinate system $\{x, y\}$
such that $\partial_x=f\tilde e_1$ and $\partial_y=k\tilde e_2$,
where $f$ and $k$ satisfy
\begin{align}
\frac{f_y}{k}=\frac{k_x}{f}, \quad \Bigl(\frac{f_y}{k}\Bigr)_y+\Bigl(\frac{k_x}{f}\Bigr)_x=fk.\label{od1}
\end{align} 
The Hamiltonian stationary condition (\ref{HS}) is equivalent to (cf. \cite{dh})
\be 
\Bigl(\frac{k}{f}\Bigr)_x+\Bigl(\frac{f}{k}\Bigr)_y=0.\label{od2}
\ee

By Theorem 3.1 of \cite{c1}, up to translations and sign,  the exact solutions of the over-determined PDE system 
 (\ref{od1})-(\ref{od2}) are given  by
\begin{align}
f=\lambda m\,{\rm csch}\Bigl(\frac{\lambda(m^2x+y)}{\sqrt{1+m^2}}\Bigr), \quad
&k=\lambda\,{\rm csch}\Bigl(\frac{\lambda(m^2x+y)}{\sqrt{1+m^2}}\Bigr); \label{cs}\\
f=\lambda m\,\sec\Bigl(\frac{\lambda(m^2x+y)}{\sqrt{1+m^2}}\Bigr), \quad
&k=\lambda\,\sec\Bigl(\frac{\lambda(m^2x+y)}{\sqrt{1+m^2}}\Bigr);\label{se}\\
f=\frac{m\sqrt{1+m^2}}{m^2x+y}, \quad &k=\frac{\sqrt{1+m^2}}{m^2x+y},\label{fra}
\end{align}
where $\lambda$ and $m$ are positive real numbers.
We note that $f=mk$. By the first equation of  (\ref{abcfk}) we have 
\be 
a=\frac{\sqrt{1+m^2}}{mk}.\nonumber
\ee
Furthermore, it follows from (\ref{change}) that 
\begin{align*}
\partial_u=\frac{1}{1+m^2}(\partial_x+m^2\partial_y),\quad
\partial_v=\frac{m}{1+m^2}(\partial_x-\partial_y).
\end{align*}

Since $a_u$ and $a_v$ are constant, the solutions (\ref{cs}) and (\ref{se}) are excluded.
 Considering also that  $a_u^2+a_v^2=1$ and $f\ne k$, we see that
 $f$ and $g$ are given by (\ref{fra}) with 
$m\ne 1$.
 In \cite[Section 7]{coz}, it is proved that the corresponding surface is congruent to 
 the Lagrangian surface obtained from case (2) with $m\ne 1$. \qed


\begin{thebibliography}{99}
\bibitem{c2} B. Y. Chen,  {\it Slant immersions,} Bull. Austral. Math. Soc. {\bf 41} (1990), 135-147.

\bibitem{c1} B. Y. Chen, {\it Solutions to over-determined systems of partial differential
equations related to Hamiltonian stationary Lagrangian surfaces,} Electron. J. Differential
Equations {\bf 2012} No. 83,  7 pp.

\bibitem{cd} B. Y. Chen, F. Dillen, {\it Warped product decompositions of real space forms
and Hamiltonian-stationary Lagrangian submanifolds,} Nonlinear Anal. {\bf 69} (2008),
3462-3494.

\bibitem{cdv} B. Y. Chen, F. Dillen and J. Van der Veken, {\it Complete classification of parallel Lorentzian surfaces in Lorentzian
complex space forms,} Internat. J. Math. {\bf 21} (2010), 665-686. 

\bibitem{cdvv} B. Y. Chen, F. Dillen, L. Verstraelen and L. Vrancken, {\it Lagrangian isometric immersions of a
real-space-form $M^n(c)$ in to a complex-space-form $\tilde M^n(4c)$,} Math. Proc. Camb. Phil. Soc. {\bf 124}
 (1998), 107-125.




\bibitem{coz} B. Y. Chen, O. J. Garay, Z. Zhou, {\it Hamiltonian-stationary Lagrangian surfaces of constant curvature $\epsilon$ in complex space forms $\tilde M^2(4\epsilon)$,} Nonlinear Anal. {\bf 71} (2009), 2640-2659.

\bibitem{co} B. Y. Chen, K. Ogiue, {\it On totally real submanifolds,} Trans. Amer. Math. Soc. {\bf 193} (1974), 257-266. 

\bibitem{dh} Y. Dong, Y. Han, {\it Some explicit examples of Hamiltonian minimal Lagrangian submanifolds in complex space forms,} Nonlinear Anal. {\bf 66} (2007), 1091-1099.


\bibitem{oh} Y. G. Oh, {\it Volume minimization of Lagrangian submanifolds under 
Hamiltonian deformations,} Math. Z. {\bf 212} (1993), 175-192.




\end{thebibliography}
\end{document}